\documentclass{amsart}
\usepackage{amssymb, amsmath}
\usepackage{latexsym}
\usepackage{graphicx}

\newcommand{\lo}{\longrightarrow}
\newcommand{\BD}{\Bbb D}
\newcommand{\BN}{\Bbb N}
\newcommand{\BC}{\Bbb C}

\newcommand{\supp}{\operatorname{supp}}

\begin{document}

\title[fourier algebra]{Fourier Algebras on Topological Foundation $*$-Semigroups}
\date{}
\author[M. Amini, A.R. Medghalchi]{Massoud Amini, Alireza Medghalchi}
\address{Department of Mathematics\\Shahid Beheshti University\\ Evin, Tehran 19839, Iran\\ m-amini@cc.sbu.ac.ir
\linebreak 
\indent Department of Mathematics and Statistics\\ University of Saskatchewan
\\ 106 Wiggins Road\\ Saskatoon, Saskatchewan\\ 
Canada S7N 5E6\\ mamini@math.usask.ca
\linebreak
\linebreak
\indent Department of Mathematics\\ Teacher Training University, Tehran, Iran}
\thanks{This research was supported by Grant 510-2090 
of Shahid Behshti University}

\begin{abstract}
We introduce the notion of the Fourier and Fouier-Stieltjes algebra of a topological $*$-semigroup and show that these are commutative Banach algebras. For a class of foundation semigroups, we show that these are preduals of von Neumann algebras. 
\end{abstract}

\maketitle

\vspace{0.25cm}
\section{Definitions and Notations} Let $S$ be a locally compact
topological semigroup and $M(S)$ be the Banach algebra of all 
bounded regular Borel measures $\mu$ on $S$ with norm and convolution
product
$$\|\mu\|=\int_S d|\mu|, \quad \int_S f(x)d(\mu \ast v)(x)=\int_S\int_S 
f(xy)d\mu(x) dv(y), (f\in C_0(S) \mu,v\in M(S)).$$
We consider the mappings $L_\mu$ and $R_\mu$ of $S$ to $M(S)$ defined by
$$L_\mu(x)=\mu \ast \delta_x , \quad R_\mu(x)=\delta_x \ast \mu \quad
(x\in S,\mu\in M(S))$$
where $\delta_x$ is the point mass at $x.$ Then the {\bf semigroup algebra}
$L(S)$ consists of those $\mu\in M(S)$ for which $L_{|\mu|}$ and
$R_{|\mu|}$ are continuous with respect to the
 weak topology of $M(S)$. $L(S)$  is
a Banach subalgebra of $M(S).$ $S$ is called {\bf foundation} if
$\cup\{\supp (\mu):\mu \in L(S)\}$ is dense in $S$.

A {\bf representation} of $S$ is a pair $\{\pi,H_\pi\}$ of a Hilbert
space $H_\pi$ and a semigroup homomorphism $\pi:S\rightarrow B(H_\pi).$
We say that $\pi$ is ({\bf weakly}) {\bf continuous} if the mappings
$x\mapsto <\pi(x)\xi,\eta>$ are continuous on $S$, for all
$\xi,\eta\in H_\pi,$ and that $\pi$ is {\bf bounded} if
$\|\pi\|={\rm sup}_{x\in S} \|\pi(x)\|<\infty.$ A {\bf topological} $*$-
{\bf semigroup} is a topological semigroup $S$ with a continuous
involution $*:S\rightarrow S.$ In this case, a representation $\pi$
is called a $*$-{\bf representation} if moreover $\pi(x^*)=\pi(x)^*
(x\in S),$ where the right hand side is the adjoint operator.

Let $S$ be a topological $*$-semigroup, $\sum(S)$ be the family of all
continuous $*$-representations $\pi$ of $S$ with $\|\pi\|\leq 1,$ and
$\sum_0(S)$ be the irreducible elements of $\sum(S).$ The $*$-representation
$\pi_u=\oplus_{\pi\in \sum_0} \pi$ in the Hilbert space $H_u=\oplus_{\pi
\in \sum_0} H_\pi$ is called the {\bf universal representation} of $S.$
This is a continuous $*$-representation with $\|\pi_u\|\leq 1.$

Now let $S$ be a foundation topological $*$-semigroup with identity.
A $*$-representation $\{\sigma, H\}$ of $L(S)$ is called {\bf non vanishing}
if for every $0\neq \xi \in H$, there exists $\mu\in L(S)$ with
$\sigma(\mu)\xi \neq 0.$ let $\sum(L(S))$ be the family of all
$*$-representations of $L(S)$ on a Hilbert space which are non vanishing,
then by theorem 2.4 and proof of corollary 2.5 in [12], one has a
bijective correspondence between $\sum(S)$ and $\sum(L(S))$ via
$$<\tilde{\pi}(\mu)\xi, \eta>=\int_S <\pi(x)\xi,\eta> d\mu(x)
\quad(\mu\in L(S),\; \xi,\eta\in H_\pi=H_{\tilde{\pi}}).$$
Moreover $\pi$ and $\tilde{\pi}$ have the same closed invariant
subspaces and so $\pi$ is irreducible if and only if $\tilde{\pi}$ is so, and
\begin{eqnarray*}
\pi(x)\tilde{\pi}(\mu)&=&\tilde{\pi}(\delta_x \ast \mu) \quad \quad
(x\in S, \mu \in L(S)).\\
\tilde{\pi}(\mu) \pi(x) &=& \tilde{\pi} (\mu \ast \delta_x)
\end{eqnarray*}
Next, following [14], we consider the family $\Omega(S)$ of all
continuous $*$-representations $\omega$ of $S$ in a
$W^*$-algebra with $\|\omega\|\leq 1$ (i.e. semigroup homomorphisms from
$S$ into the unit ball of a $W^*$-algebra $M_\omega$ such that
$\omega(x^*)=\omega(x)^*$ for all $x\in S$).

\vspace{0.25cm}
\section{Semigroup $C^*$-algebra and big semigroup algebra}
Let $S$ be a topological $*$-simgroup with identity. Given $\rho\subseteq
\sum=\sum(S)$ and $\mu\in L(S),$ define $\|\mu\|_\rho={\rm sup} \{
\|\tilde{\pi}(\mu)\|:\pi \in \rho \}$ and $I_\rho=\{\mu\in L(S):
\|\mu\|_\rho=0\}.$ Then $I_\rho$ is clearly a closed two-sided ideal of
$L(S)$ and $\|\mu+I_\rho\|=\|\mu\|_\rho$ defines a $C^*$-norm on
$L(S)/I_\rho.$ The completion of this quotient space in this norm is
a $C^*$-algebra which is denoted by $C_\rho ^*(S).$  When
$\rho=\sum,$ then $C^*$-algebra $C^*(S)=C^*_\Sigma(S)$ is called
the ({\bf full}) {\bf semigroup $C^*$-algebra} of $S$.
If $S$ is foundation and $\sum$ separates the points of $S$, then $L(S)$
is $*$-semisimple [12], and so $I_\Sigma=\{0\}.$ In this case $L(S)$ is a
norm dense subalgebra of $C^*(S).$ In particular each $\tilde{\pi}\in \sum
(L(S))$ uniquely extended to a $*$-representation of $C^*(S)$ which we
still denote it by $\tilde{\pi}.$

Following [14] one may also consider the $C^*$-algebra $C^*_\Omega (S)$
generated by the set of operators $\{\pi(x):x\in S, \pi \in \Omega(S)\}.$
If $\omega_\Omega= \oplus_{\omega\in\Omega} \omega$ and $M_\Omega= \sum_
{\omega\in \Omega} \oplus M_\omega$, then $C^*_\Omega(S)$ is the span
closure of $\omega_\Omega(S)$ in $M_\Omega.$ We shall define the
Fourier-Stieltjes algebra $B(S)$ naturally as the linear span of all
continuous positive definite functions on $S$ and one expects that
as in the group case [8], the Fourier-Stieltjes algebra is isometically
isomorphic to the dual space of the semigroup $C^*$-algebra. But
$C^*_\Omega (S)^*$ only contains an isometic copy of $B(S)$ in general
(see remark 4.3.(e) in [13]), and $C^*_\Omega (S)^*\simeq B(S),$
when $S$ is discrete. As we would show in the next section, if one replaces
$C^*_\Omega(S)$ with $C^*(S)$, this pathology disapears and one gets
$C^*(S)^*\simeq B(S)$ for each foundation topological *-semigroup
with identity for which $\sum(S)$ separates the points of $S$ (theorem 3.3).

{\bf Warning:} The $C^*$-algebra $C^*(S)$ defined above is not the same
as the one defined in [5], [6]. The reason is that in the latter, one
consideres representations instead of $*$-representation. Doing so,
$C^*(S)$ could be non commutative even for abelian discrete semigroups,
e.g.$ C^*(\BN)$ turns out to be the Toeplitz algebra [16],
where as here $C^*(\BN)$ (when $\BN$ is endowed with trivial
involution $n^*=n$) is nothing but the algebra $C(\BD)$ of
continuous functions on the unit disk.

Let $W^*(S)$ be the enveloping $W^*$-algebra of $C^*(S).$ If
$\{\pi_u,H_u\}$ is the universal $*$-representation of $S$ then
$\{\tilde{\pi}_u,H_u\}$ is the universal representation of 
$C^*(S)$ and $W^*(S)=\tilde{\pi}_u(C^*(S))''\subseteq B(H_u)$ is isometrically 
isomorphic to $C^*(S)^{**}$ with Arens product [17]. On the other hand, let 
$W^*_\Omega(S)=<\omega_\Omega(S)>^{-\sigma}$
 be the weak operator closure of the 
linear span of $\omega_\Omega(S)$
 in $M_\Omega$. Unlike the case of $C^*$-algebras, 
one can easily see that the $W^*$-algebras $W^*(S)$ and $W^*_\Omega(S)$ 
coincide in general. Indeed to each $\pi\in\sum(S)$, there corresponds a 
*-representation $\omega_\pi:S\lo M_\pi=<\pi(S)>^{-\sigma}\subseteq B(H_\pi)$ 
given by $\omega_\pi(x)=\pi(x)$ $(x\in S)$. Conversely, if $\omega\in\Omega(S)$ 
then one can take the universal representation $\{i_\omega,H_\omega\}$ of 
$M_\omega$ to embed $M_\omega$ 
isometrically in $B(H_\omega)$. Then $\pi_\omega=i_\omega\circ
 \omega\in\sum(S)$, 
$M_{\pi_\omega}=<\pi(S)>^{-\sigma}=<i_\omega(\omega(S)>^{-\sigma}=M_\omega$, 
and $\omega_{\pi_\omega}(x)=i_\omega(x)=\omega(x)$ $(x\in S)$. 

In the next proposition,
 by an abuse of notation we use $\sum(S)$ to denote the 
set of operators $\{\pi(s):\pi\in\sum, s\in S\}$ (rather than the set of 
representations). Also $\sum(L(S))=\{\tilde{\pi}(\mu): \pi\in\sum, \mu\in 
L(S)\}$ and $\sum(C^*(S))$ is defined similarly. 

\vspace*{0.5cm}
\noindent
{\bf Proposition 2.1.} If $S$ is a foundation topological *-semigroup with 
identity and $\sum$ separates the points of $S$, then
$$\sum(S)''=\sum(L(S))''=\sum(C^*(S))''=W^*(S),$$
where the commutators are taken in the corresponding universal representations.

{\bf Proof.} Since $\sum$ separates the points of $S$, $L(S)$ is norm dense in 
$C^*(S)$ and so $\sum(L(S))$ is ultra-weakly dense in $\sum(C^*(S))$, hence 
$\sum(L(S))''=\sum(C^*(S))''=W^*(S)$. Now
$$<\tilde{\pi}(\mu)\xi,\eta>=\int_S<\pi(x)\xi,\eta>d\mu(x)\quad (\mu\in 
L(S),\xi,\eta\in H_\pi, \pi\in \sum).$$
Hence if $T\in B(H_u)$ commutes with all operators $\pi(x)$, then it also 
commutes with all operators $\tilde{\pi}(\mu)$, that is $\sum(S)'\subseteq 
\sum(L(S))'$. Hence $\sum(S)''\supseteq\sum(L(S))''=W^*(S)$. On the 
other hand, for each $x\in S$,  
$\tilde{\pi}(\delta_x)=\pi(x)$, so 
$\sum(S)''\subseteq\sum(\tilde{L}(S))''=W^*(S)$ and the equality holds. 
Note that we have used the same notation
 $\tilde{\pi}$ for the extension of $\tilde{\pi}$ to 
$M(S)$.$\Box$

\vspace*{0.5cm}
\noindent
{\bf Proposition 2.2.} Let $S$ be a foundation topological *-semigroup with 
identity and $\sum(S)$ separates the points of $S$, then $L(S), M(S)$
and $C^*(S)$ are isometrically embedded in $W^*(S)$.

{\bf Proof.} Consider the universal representation $\{\pi_u,H_u\}$ of $S$. 
Since $\sum(S)$ separates the points of $S$, this is a faithful 
representation. Hence by theorem 2.4 of [12],
 the correseponding representation 
$\tilde{\pi}_u$ of $L(S)$ is faithful, which in turn is extended to a 
faithful representation of $C^*(S)$. Also $M(S)$ is the multiplier algebra of 
$L(S)$ and $L(S)$ has a bounded approximate identity, so 
$\tilde{\pi}_u$ is extended
 to a faithful representation of $M(S)$. The theorem 
follows now from Proposition 2.1. $\Box$

Next we study the $C^*$-algebras $C^*_\rho(S)$, for any $\rho\subseteq\sum$ 
and foundation semigroup $S$. Consider $\pi\in\sum(S)$, then $\tilde{\pi}$ is
extended to a *-representation of $M(S)$ on $H_\pi$ with 
$$\tilde{\pi}(\mu)=\int_S\pi(x)d\mu(x)\quad (\mu\in M(S)).$$
If $\rho\subseteq\sum$, then 
$\|\mu\|_\rho=\sup_{\pi\in\rho}\|\tilde{\pi}(\mu)\|$ defines a seminorm on 
$M(S)$ and for each $\mu,\nu\in M(S)$,
$$\|\mu*\nu\|_\rho\leq\|\mu\|_\rho\|\nu\|_\rho,\; \|\mu\|_\rho=\|\mu^*\|_\rho, 
\;\|\mu*\mu^*\|_\rho=\|\mu\|_\rho^2,\; \|\mu\|_\rho\leq\|\mu\|_1.$$
Let $N_\rho=\{\mu\in L(S): \|\mu\|_\rho=0\}$, then for 
$\dot{\mu}=\mu+N_\rho\in L(S)/N_\rho$, 
$\|\dot{\mu}\|_\rho=\|\mu\|_\rho$defines a norm on $L(S)/N_\rho$ and is 
independent of the representative $\mu$ of the class $\dot{\mu}$, and the 
completion of $L(S)/N_\rho$ is the $C^*$-algebra $C_\rho^*(S)$. The next 
result and its proof are similar to the corresponding result in [8].

\vspace*{0.5cm}
\noindent
{\bf Proposition 2.3.} For $\rho\subseteq\sum$, let 
$N'_\rho=\bigcap_{\pi\in\rho}\ker\tilde{\pi}\subseteq C^*(S)$, then 
$N'_\rho\supseteq N_\rho$ and the map $\mu\mapsto\dot{\mu}$ of $L(S)$ onto 
$L(S)/N_\rho$ is extended uniquly to an isometric epimorphism: $C^*(S)\lo 
C^*_\rho(S)$ with kernel $N'_\rho$; i.e. $C^*_\rho(S)\simeq C^*(S)/N'_\rho$.

{\bf Proof.} Since $\|\dot{\mu}\|_\rho\leq\|\mu\|_\Sigma$, the map 
$\mu\mapsto\dot{\mu}$ is extended continuously to an *-homomorphism 
$\pi_\rho:C^*(S)\lo C^*_\rho(S)$ which decreases the norm. Since the image of 
$\pi_\rho$ contains $L(S)/N_\rho$, it is dense in $C^*_\rho(S)$. On 
the other hand, the range of $C^*$-algebra homomorphisms are always closed 
[17], so $\pi_\rho$ is surjective. If $g\in N'_\rho$, then there is a net 
$\{\mu_\alpha\}\subseteq L(S)$ such that $\|g-\mu_\alpha\|_\Sigma\lo 0$. 
Then $\|\mu_\alpha\|_\rho=\sup_{\pi\in\rho}\|\tilde{\pi}(\mu_\alpha
)\|=\sup_{\pi\in\rho}\|\tilde{\pi}(g-\mu_\alpha)\|\leq\|g-
\mu_\alpha\|_\Sigma\lo 0$, hence 
$\pi_\rho(g)=\lim_\alpha\pi_\rho(\mu_\alpha)=0$, i.e. 
$N'_\rho\subseteq\ker\pi_\rho$. Conversely, if $g\in C^*(S)$ and 
$\pi_\rho(g)=0$, then there is a net $\{\mu_\alpha\}\subseteq L(S)$
such that $\|\dot{\mu}_\alpha\|_\rho=\|\mu_\alpha\|_\rho\lo 0$ and 
$\|g-\mu_\alpha\|_\Sigma\lo 0$, hence 
$\sup_{\pi\in\rho}\|\tilde{\pi}(g)\|=\lim_\alpha(\sup_{\pi\in\rho}\| 
\tilde{\pi}(g)-\tilde{\pi}(\mu_\alpha)\|)\leq\lim_\alpha\|g-
\mu_\alpha\|_\Sigma=0$, i.e. $\tilde{\pi}(g)=0$ for all $\pi\in\rho$, that's 
$g\in N'_\rho$. Therefore $\ker \pi_\rho=N'_\rho$. Finally  
$\pi_\rho$ is clearly surjective. $\Box$

{\bf Remark.} $\rho\subseteq\sum$ separates the points of $S$ exactly when 
$N_\rho=0$. This condition is generally weaker than $N'_\rho=0$. 

\vspace*{0.5cm}
\section{Fourier and Fourier Stieltjes algebra}
Let $S$ be a unital topological *-semigroup. A complex valued function 
$u:S\lo\BC$ is called {\bf positive definite} if for all positive integers $n$ 
and all $\lambda_1,\dots,\lambda_n\in\BC$, and $x_1,x_2,\dots,x_n\in S$, we have
$$\sum_{i=1}^n \sum_{j=1}^n \lambda_i\bar{\lambda}_j u(x_ix_j^*)\geq 0.$$
Let $P(S)$ denotes the set of all continuous positive definite functions on 
$S$. We denote the linear span of $P(S)$ by $B(S)$ and call it the {\bf 
Fourier Stieltjes algebra} of $S$.

\vspace*{0.5cm}
{\bf Proposition 3.1.} Let $S$ be a unital topological *-semigroup and 
$u:S\lo\BC$ be a bounded function, then the followings are equivalent:
\begin{itemize}
\item[(i)] $u\in P(S)$
\item[(ii)] There is a $\sigma$-continuous *-representation $\pi:S\lo B(H)$ 
and a cyclic vector $\xi\in H$ such that $u=<\pi(\cdot)\xi,\xi>$ 
and $\|\pi\|\leq 1$.
\item[(iii)] There is a $W^*$-algebra $M$ and a $\sigma$-continuous 
*-representation $\omega:S\lo
 M$ into the unit ball of $M$, and a positive bounded 
normal linear functional $\psi\in M_*$ such that $u=\psi\circ \omega$ and 
$<\omega(S)>^{-\sigma}=M$. 
\end{itemize}

{\bf Proof.} The proof of equivalence of (i) and (iii) is contained in [14]. 
Also (ii)$\Longrightarrow$(i) is trivial. If (iii) holds, then by 
Gelfand-Naimark-Segal-construction,
there corresponds to $\psi\in M_{*+}$, a nondegenerate
$\sigma$-continuous *-representation $\{\pi_\psi,H_\psi,\xi_\psi\}$ of $M$ 
such that $\psi=<\pi_\psi(\cdot)\xi_\psi,\xi_\psi>$, and 
$\overline{\pi_\psi(M)\xi_\psi}=H_\psi$. Let $\pi=\pi_\psi\circ\omega$,
then $\pi$ 
is a $\sigma$-continuous *-representation of $S$ in $H_\psi$, and 
$u=\psi\circ \omega=<\pi_\psi\circ \omega(\cdot)\xi_\psi,\xi_\psi>= 
<\pi(\cdot)\xi_\psi,\xi_\psi>$, and 
$\pi(S)\xi_\psi=\pi_\psi(\omega(S))\xi_\psi$
is $\sigma$-total in $H_\psi$. Finally 
$\|\pi(x)\|=\|\pi_\psi(\omega(x))\|\leq\|\omega(x)\|\leq
1$, for each $x\in S$, and (ii) follows. $\Box$

Now consider the algebra $F(S)$ defined by Lau in [14].

\vspace*{0.5cm}
{\bf Corollary 3.2.} If $S$ is a unital topological *-semigroup, then 
$B(S)=F(S)$. $\Box$

\vspace*{0.5cm}
{\bf Theorem 3.3.} Let $S$ be a foundation topological $*$-semigroup with 
identity such that $\sum(S)$ separates the points of $S$. Then $B(S)$ and 
$C^*(S)^*$ are isometrically isomorphic as Banach spaces. 

\vspace*{0.5cm}
{\bf Proof.} Since $\sum(S)$ separates the points of $S$, $L(S)$ is 
$*$-semisimple [11]. In particular $I_\Sigma=\{0\}$, and so $L(S)\subseteq 
C^*(S)$ is norm dense. Therefore each continuous $*$-representation
 of $L(S)$ uniquely 
extends to one on $C^*(S)$. Now each $u\in B(S)\subseteq C_b(S)$ acts on 
$L(S)$ via 
$$<u,\mu>=\int_S u(x)d\mu(x)\quad (\mu\in L(S)).$$
On the other hand, to each $u\in B(S)$, there corresponds some 
$\pi\in\Sigma(S)$ and $\xi,\eta\in H_\pi$ such that $u(x)=<\pi(x)\xi,\eta>$ 
(Proposition 3.1), and 
\begin{align*}
\sup\{|<u,\mu>|: \mu\in L(S), \|\mu\|_\Sigma\leq 1\}& = 
\sup_{\|\mu\|_\Sigma\leq 1}|\int_S u(x)d\mu(x)|\\
&=\sup_{\|\mu\|_\Sigma\leq 1}|\int_S<\pi(x)\xi,\eta>d\mu(x)|\\
&=\sup_{\|\mu\|_\Sigma\leq 1}|<\tilde{\pi}(\mu)\xi,\eta>|\\
&\leq\|\tilde{\pi}\| \|\xi\| \|\eta\| \leq \|\xi\|.\|\eta\|.
\end{align*}
Hence $u$ extends uniquely to an element of $C^*(S)^*$. Also by density of 
$L(S)$ in $C^*(S)$, the norm of $u$ as an element of $B(S)$ is the same as the 
norm of $u$ as an element of $C^*(S)^*$. 

Now let $T\in C^*(S)^*_+$, then by restriction, this gives an element of 
$L(S)^*_+$, hence by the proof of theorem 2.4 (for $\omega=1$) of [11], 
there is $\pi\in\Sigma(S)$ and $\xi\in H_\pi$ such that 
$$T(\mu)=\int_S<\pi(x)\xi,\xi>d\mu(x)\quad (\mu\in L(S)).$$
Define $u\in C_b(S)$ by $u(x)=<\pi(x)\xi,\xi>$, then $u\in P(S)$ and 
$<u,\mu>=\int_S<\pi(x)\xi,\xi>d\mu(x)=<T,\mu>$ $(\mu\in L(S))$, and so $T=u\in 
P(S)$. $\Box$  

\vspace*{0.5cm}
{\bf Theorem 3.4.} If $S$ is a unital foundation topological *-semigroup such 
that $\sum(S)$ separates the points of $S$, then the following norms coincides 
on $B(S)$ and $B(S)$ is a commutative Banach algebra under each of these norms 
and $B(S)\simeq W^*(S)_*$, isometrically isomorphic as Banach spaces:
\begin{itemize}
\item[(i)] $\|u\|=\sup\{|\int_S u(x)d\mu(s)| :\mu\in L(S), 
\sup_{\pi\in\sum} \|\tilde{\pi}(\mu)\|\leq 1\}$,
\item[(ii)] $\|u\|=\inf\{\| \psi\|: \psi\in M_* \; \text{and}\; 
u=\psi\circ\omega,\; \text{for some}\; (\omega,M)\in\Omega(S)\},$
\item[(iii)] $\|u\|=\sup\{|\sum_{n=1}^N \lambda_n u(x_n)|: N\geq 1, 
\lambda_n\in\BC, x_n\in S\; (1\leq n\leq N)$, \\
$\sup_{\pi\in\sum}\|\sum_{n=1}^N\lambda_n\pi(x_n)\|\leq 1\}$,
\item[(iv)] $\|u\|=\inf\{\|\xi\|\|\eta\|: u=<\pi(\cdot)\xi,\eta>,\;\text{for 
some}\; \{\pi,H_\pi\}\in\sum, \xi,\eta\in H_\pi\}.$
\end{itemize}

{\bf Proof.} It is shown in [14] that $F(S)\simeq W^*_\Omega(S)_*$, with 
respect to the 
norm (ii) and in Proposition 3.3 that $B(S)\simeq W^*(S)_*$, with respect to
 the norm (i). 
Now $F(S)=B(S)$ and $W^*_\Omega(S)$ and $W^*(S)$ are isometrically isomorphic 
as $W^*$-algebras. Hence, by the uniqueness of the predual, the norms (i) and 
(ii) are the same. Next, by Proposition 2.1, $W^*(S)$ is the $W^*$-algebra 
generated by $\{\tilde{\pi}(\delta_x):x\in S, \pi\in\sum(S)\}$. Hence the 
linear combinations of the form $\sum_{n=1}^N\lambda_n\delta_{x_n}$ are 
$\sigma$-dense in $W^*(S)$. Hence, by Kaplansky's density theorem [17], 
norms (i) and (iii) coincide. Finally, if $u=<\pi(\cdot)\xi,\eta>$, for some 
$\pi\in\sum(S)$ and $\xi,\eta\in H_\pi$, then 
\begin{align*}
\|u\|&=\sup_{\|\mu\|_\Sigma\leq 1}|\int_S ud\mu|=\sup_{\|\mu\|_\Sigma\leq 1} 
|\int_S<\pi(x)\xi,\eta>d\mu(x)|=\sup_{\|\mu\|_\Sigma\leq 1}| 
<\tilde{\pi}(\mu)\xi,\eta>| \\
 &\leq\sup_{\|\mu\|_\Sigma\leq 1} \|\tilde{\pi}(\mu)\|\|\xi\|\|\eta\| \leq\| 
\xi\|\|\eta\|.
\end{align*}

On the other hand, given $\varepsilon>0$, there is $(\omega,M)\in\Omega(S)$
 and 
$\psi\in M_*$ such that $u=\psi\circ w$ and $\|\psi\|\leq\|u\|+\varepsilon$. 
Now, by Gelfand-Naimark-Segal-construction,
 there is a representation $\{\pi_\psi, H_\psi\}$ of 
$M$ and vector $\xi,\eta\in H_\psi$ such that 
$\psi=<\pi_\psi(\cdot)\xi,\eta>$ and $\|\psi\|=\|\xi\|\|\eta\|$. Put 
$\pi=\pi_\psi\circ\omega$, then $\{\pi,H_\psi\}\in \sum(S)$ and $u=\psi\circ 
\omega=<\pi_\psi\circ\omega(\cdot)\xi,\eta>=<\pi(\cdot)\xi,\eta>$,
 and $\|u\|\geq 
\|\psi\|-\varepsilon=\|\xi\|\|\eta\|-\varepsilon$. Hence norms (ii) and (iv) 
also coincide. Now by theorems 3.2 and 4.1 of [14], $B(S)$ is a Banach algebra 
under pointwise multiplication and norm (ii) and by Theorem 3.3, 
$B(S)\simeq W^*(S)_*$, under the norm (i), and we are done. $\Box$

\vspace*{0.5cm}
{\bf Proposition 3.5.} Let $S$ be a foundation topological *-semigroup with 
identity, then the followings are equivalent.
\begin{itemize}
\item[(i)] $\sum(S)$ separates the points of $S$,
\item[(ii)] $B(S)$ separates the points of $S$,
\item[(iii)] $L(S)$ is *- semisimple,
\item[(iv)] $M(S)$ is *-semisimple,
\end{itemize}

{\bf Proof.} This follows from our Proposition 3.1, and Lemma 3.1 and Theorem 
3.4 of [12]. $\Box$ 

Let $S$ be a topological *-semigroup and 
$C_c(S)$ be the algebra of all continuous functions on $S$ of compact support. 
Then the closed subalgebra $\overline{(B(S)\cap
 C_c(S))}\subseteq B(S)$ is denoted by 
$A(S)$ and is called the {\bf Fourier algebra} of $S$. Since $B(S)\cap C_c(S)$ 
is clearly a two-sided ideal of $B(S)$ under pointwise multiplication, we have 
the following result.

\vspace*{0.5cm}
{\bf Proposition 3.6.} $A(S)$ is a closed two-sided ideal of $B(S)$, in 
particular $A(S)$ is a commutative Banach algebra under pointwise 
multiplication and each of the norms in Theorem 3.4. $\Box$

We observed that $B(S)$ is the predual of a von Neumann algebra. Now we want 
to see when this holds for $A(S)$.
The following lemma should be a part
 of the literature, but we give a proof for 
the sake of completness.

\vspace*{0.5cm}
{\bf Lemma 3.7.} Assume that $W$ is a von Neumann algebra,
$B$ is the predual Banach space of $W$, 
$z_0\in W$ is a central projection, and 
$A=z_0.B$ (Arens product).
Then $A^*\simeq z_0B^*$ (isometrically 
isomorphic as Banach spaces).

\vspace*{0.5cm}
{\bf Proof.} It is clear that $A\subseteq B$. Let $j:A\lo B$ be the 
corresponding linear inclusion, then $j^*:B^*\lo A^*$ is a surjective 
continuous linear map and $A^*\simeq B^*/\ker j^*$. Now given $a\in A$ and 
$f\in B^*$, there is $b\in B$ such that $a=z_0.b$, and so 
$(1-z_0).a=(1-z_0).(z_0.b)=((1-z_0)z_0).b=0$, hence
$$<j^*((1-z_0)f),a>=<(1-z_0)f, j(a)>=<f, (1-z_0).a>=0.$$
Hence $(1-z_0)B^*\subseteq\ker j^*$. Conversely if $f\in \ker j^*$, then for 
each $b\in B$,
$$<z_0f,b>=<f,z_0.b>=<f,j(z_0.b)>=<j^*(f),z_0.b>=0$$
hence $z_0f=0$ in $B^*$; so $f=(1-z_0)f\in(1-z_0)B^*$. Therefore $\ker 
j^*=(1-z_0)B^*$ and so $B^*/\ker j^*\simeq z_0B^*$. $\Box$

{\bf Lemma 3.8.} Let $S$ be a foundation topological $*$-semigroup with 
identity such that $\sum(S)$ separates the points of $S$, then $A(S)$ is 
uniformely dense in $C_0(S)$.

\vspace*{0.5cm}
{\bf Proof.}  Let $E=B(S)\cap C_c(S)$. If there is $x\in S$ such that 
$u(x)=0$, for each $u\in B(S)$, then for $\delta_x\in W^*(S)$ (see Proposition 
2.2) 
$$<\delta_x,u>=u(x)=0 \quad (u\in B(S))$$
i.e. $\delta_x=0$ in $W^*(S)$, which is a contradiction. Hence there is $u\in 
B(S)$ such that $u(x)\neq 0$, say $u(x)=1$. Now let $v\in E$ be such that 
$\|u-v\|<\frac{1}{2}$, then $|u(x)-v(x)|\leq \|u-v\|_\infty\leq 
\|u-v\|<\frac{1}{2}$, so $|v(x)|\geq \frac{1}{2}$. Therefore $E$ 
Vanishes nowhere on $S$. Next, 
we show that $E$ separates the points of $S$. Let $x,y\in S$ and $x\neq y$. 
Then, by Proposition 3.5, there is $u\in B(S)$ such that $u(x)\neq u(y)$. 
Choose $v\in E$ such that $v(x)\neq 0$, say again $v(x)=1$. If $v(y)\neq 1$, 
we have nothing to proof. Otherwise, $v(x)=v(y)=1$ and $(uv)(x)\neq (uv)(y)$, 
and we are done. Now the result follows from Stone-Weierstauss  theorem. 
$\Box$ 

\vspace*{0.5cm}
{\bf Theorem 3.9.} Let $S$ be a foundation locally compact 
topological *-semigroup with 
identity such  that $\sum(S)$ separates the points of $S$,
then the followings are equivalent.
\begin{itemize}
\item[(i)] $A(S)\subseteq C^*(S)^*=B(S)$ is invariant under Arens product (i.e.
$u.g, g.u\in 
A(S)$, for each $u\in A(S)$, $g\in C^*(S)$).
\item[(ii)] $A(S)\subseteq W^*(S)_*=B(S)$ is invariant under Arens product 
(i.e. $u.T, T.u\in A(S)$, for each $u\in A(S)$, $T\in W^*(S)$).
\item[(iii)] $A(S)=I^{\perp}$, for some $w^*$-closed
 ideal $I$ of $W^*(S)$.
\item[(iv)] $A(S)=z_0$. $B(S)$, for some central projection $z_0\in W^*(S)$.
\item[(v)] $A(S)$ is translation invariant. 
\item[(vi)] $A(S)^*$ is a $W^*$-algebra (indeed a quotient of $W^*(S)$).
\item[(vii)] For each compact subset $K\subseteq S$ and each $x\in S$, both 
$Kx^{-1}$ and $x^{-1}K$ are compact. 

Moreover if $S$ is not compact, then these are equivalent 
to:
\item[(viii)] The one-point compactification $S_\infty$ of $S$ is a semigroup 
compactification of $S$.
\item[(ix)] $C_0(S)$ is translation invariant. 
\item[(x)] The $C^*$-algebra $F=C_0(S)\oplus \BC$ is an $m$-admissible 
subalgebra of $S$ and $S_\infty$ is an $F$-compactification of $S$ (see 
[3]). 
\end{itemize}

\vspace*{0.5cm}
{\bf Proof.} The equivalence of (i), (ii), (iii), and (iv) holds in general 
[4]. Also equivalence of (vii), (viii), (ix), and (x) for the case of 
locally compact, noncompact semigroups is contained in [3]. To complete the 
proof, it is enough to show that $(iv)\Longrightarrow(v)\Longrightarrow(iii)$,
$(iv)\Longrightarrow(vi)\Longrightarrow(iii)$, and $(vii)\Longrightarrow(v)$:

$(iv)\Longrightarrow(v)$. Given $u\in A(S)$, there is $v\in B(S)$ such that 
$u=v.z_0$. For each $x\in S$, consider $\delta_x\in M(S)\subseteq 
W^*(S)=B(S)^*$, then 
$$_xu=u.\delta_x=(v.z_0).\delta_x=v.(\delta_xz_0)=v.(z_0\delta
_x)=(v.\delta_x).z_0=(_xv).z_0\in A(S),$$
since $B(S)$ is always left translation invariant (Theorem 3.2 of [14]).

$(v)\Longrightarrow(ii)$. By assumption, $\delta_x.u$, $u.\delta_x\in A(S)$, 
for each $u\in A(S)$, $x\in S$. But finite sums of the form $\sum_{i=1}^N 
c_i\delta_{x_i}$ are $w^*$-dense in $W^*(S)$ (Proposition 2.1), and
Arens product is $w^*$-continuous, hence we get that $A(S)$ is invariant under 
Arens product by elements of $W^*(S)$.

$(iv)\Longrightarrow(vi)$. If $A(S)=z_0.B(S)$, then $A(S)^*\simeq z_0W^*(S)$ 
is a quotient of the $W^*$-algebra $W^*(S)$ (lemma 3.7).

$(vi)\Longrightarrow(iii)$. If $A(S)^*$ is a quotient of $W^*(S)$, then there 
is a central projection $z_0\in W^*(S)$ such that $A(S)^*\simeq z_0 W^*(S)$. 
Take $I=(1-z_0)W^*(S)$, then $(I^\perp)^*\simeq W^*(S)/I\simeq z_0W^*(S)\simeq 
A(S)^*$, so by the uniquencess of the predual of a $W^*$-algebra, $A(S)\simeq 
I^\perp$ as Banach spaces.

$(vii)\Longrightarrow(v)$. Given $u\in A(S)$, there is a net 
$\{u_\alpha\}\subseteq B(S)\cap C_c(S)$ such that $u_\alpha\lo u$ in $B(S)$. 
By the $w^*$-continuity of the Arens product, for a given $x\in S$, we have 
$_xu_\alpha=\delta_x.u_\alpha\lo\delta_x.u=_x\!u$, in the norm of $B(S)$. 
But $\supp(_xu_\alpha)=x^{-1}(\supp u_\alpha)$ is compact by assumption, 
that's $_xu_\alpha\in B(S)\cap C_c(S)$ for each $\alpha$, hence $_xu\in A(S)$ 
as required.

$(v)\Longrightarrow(vii)$ 
Let $E=B(S)\cap C_c(S)$, then $E$ is dense in $A(S)$ by definition. Also, by 
the proof of lemma 3.8, $E$ is dense in $C_0(S)$ with respect to the uniform 
norm. Now given compact subset $K\subseteq S$, there is $g\in C_0(S)$  
such that $g=1$ on $K$ (Uryshon's lemma). Take $h\in E$ with $\|h-g\|_\infty
<\frac{1}{2}$.
 Then clearly $|h|\geq \frac{1}{2}$ on $K$. By assumption, $_xh\in 
A(S)$, and so there is $k\in E$ such that in the norm of $A(S)$ we have 
$\|k-_xh\|<\frac{1}{4}$. But $\|k-_xh\|_\infty\leq\|k-_xh\|<\frac{1}{4}$ (by 
Theorem 3.2(b) of [4]), and $|_xh|\geq\frac{1}{2}$ on $x^{-1}K$, hence 
$|k|\geq\frac{1}{4}$ on $x^{-1}K$. In particular $x^{-1}K\subseteq\supp(k)$. 
Since $x^{-1}K$ is compact, therefore $x^{-1}K$ is compact.$\Box$

We call a topological semigroup $S$ ({\bf left}) {\bf weakly concellative 
} if it satisfies (the first condition of (iii))
 condition (iii) of the above 
theorem. Note that if $S$ is a Hausdorff cancellative topological semigroup, 
then it is weakly concellative (the 
injective continuous map $y\mapsto yx$ would be a homeomorphism of the compact 
set $K$ onto $Kx^{-1}$).

If $S$ is a locally compact, left weakly concellative, Hausdorff topological 
semigroup, then it is well known that $S$ has a right Haar measure if and 
only if it 
contains a unique minimal left ideal (which is necessarily closed) [1].

If $S$ is a unital, foundation, left weakly concellative topological 
*-semigroup, then by the above proposition, $A(S)^*$ is a $W^*$-algebra 
quotient of $W^*(S)$. We denote this $W^*$-algebra by $VN(S)$ and call it the 
{\bf semigroup von Neumann algebra} of $S$. 

\vspace*{0.5cm}
{\bf Corollary 3.9.} If $S$ is as above, then $W^*(S)\simeq A(S)^\perp\oplus 
VN(S)$; as $W^*$-algebras.

\vspace*{0.5cm}
{\bf Proof.} We have
$$A(S)^\perp\simeq (B(S)/A(S))^*\simeq (B(S)/z_0.B(S))^*\simeq 
(1-z_0)B(S)^*\simeq(1-z_0)W^*(S),$$
and $VN(S)\simeq z_0W^*(S)$.$\Box$

Next we consider the problem that when $B(S)$ is the multiplier algebra of 
$A(S)$. (This is always true for amenable
 topological groups.) First some notations: 
Let $A$ be a Banach algebra and consider $A^{**}$ with first Arens
 product. This is 
also a Banach algebra which contains (a copy of) $A$. Then
\begin{align*}
M(A) & =\{x\in A^{**}: xA\cup Ax\subseteq A\}\\
\tilde{M}(A)&=\{x\in (A^*.A)^*: xA\cup Ax\subseteq A\}
\end{align*}
are Banach subalgebras of $A^{**}$ containing $A$. $M(A)$ is called the 
multiplier algebra of $A$. If $A$ has a left bounded approximate identity, 
then $\tilde{M}(A)=RM(A)$, where $R$ stands for right and $RM(A)=\{x\in 
A^{**}:Ax\subseteq A\}$. In particular, if $A$ moreover is commutative, 
then $\tilde{M}(A)=M(A)$.

The proof of the next theorem is a modification of the argument 
given in [15] for topological groups,
 using the central projection $z_0$ of the 
above proposition. 

\vspace*{0.5cm}
{\bf Theorem 3.10.} Let $S$ be a unital, foundation, left weakly concellative 
topological *-semigroup such that $\sum(S)$ separates the points of $S$. Then 
if $A(S)$ has a bounded approximate identity, we have $B(S)=M(A(S))$.

\vspace*{0.5cm}
{\bf Proof.} Consider the map $T:B(S)\lo (A(S)^*.A(S))^*$ defined by 
$$<T(a), F.a>=<F, ab> \quad (a\in B(S), b\in A(S), F\in A(S)^*).$$
Then $T$ is clearly a linear isometry, and so one may consider $B(S)$ as a 
Banach subspace of $(A(S)^*.A(S))^*$,
 where the action on $A(S)$ is the restriction 
of the product of $B(S)$. Since $A(S)$ is a two sided (closed) ideal of 
$B(S)$, we have $B(S)\subseteq\tilde{M}(A(S))$. Conversely, if 
$h\in\tilde{M}(A(S))$, then $h$ is a linear operator on the closed subspace 
$A(S)^*.A(S)$ of $A(S)^*$, and so it is
 extended to a linear operator $\tilde{h}$ 
(with the same norm) on $A(S)^*$ (by Hahn-Banach theorem) such that 
$a.h=a.\tilde{h}\in A(S)$, for each $a\in A(S)$. Let $\{e_\alpha\}\subseteq 
A(S)$ be a bounded approximate identity, then passing to a subset (if 
necessary) we may assume that $e_\alpha\overset{w^*}{\lo} e$, the identity of 
$A(S)^{**}$, and that $e_\alpha.\tilde{h}\overset{w^*}{\lo} \lambda$, for some 
$\lambda\in B(S)=C^*(S)^*$. Then, given $g\in C^*(S)$, we have 
$$<g,\lambda>=\lim_{\alpha}<g,e_\alpha.h>=\lim_\alpha<g,e_\alpha.\tilde{h} 
>=<g,e.\tilde{h}>.$$
But, given $a\in A(S)$, we have
$$<g.a,e.\tilde{h}>=<g,a(e.\tilde{h})>= <g,(ae).\tilde{h}> =<g,a.\tilde{h}> = 
<g,a.h>=<g.a,h>,$$
hence $e.\tilde{h}=h$, and so $<g,\lambda>=<g,h>$, for each $g\in C^*(S)$.

Next, given $g\in C^*(S)$ and $a\in A(S)$, we have $g.a\in C^*(S)$ (note that 
we are now talking about the Arens triple $B(S)$, $B(S)^*$, and $B(S)^{**}$, 
where as in 3.8(i) we were talking about the triple $C^*(S)$, $C^*(S)^*$, 
and $C^*(S)^{**}$). Indeed, this follows from the fact that the product of 
$B(S)$ is a separately $w^*$-continuous bilinear map. Hence
$$<g,a\lambda>=<g.a,\lambda>=<g.a,h>=<g,ah> \quad (g\in C^*(S), a\in A(S)).$$
But $a\lambda,ah\in A(S)$ (since $A(S)$ is an ideal of $B(S)$ and $h$ is a 
multiplier of $A(S)$), and $A(S)=z_0.B(S)$, hence $z_0.a\lambda=a\lambda$, and 
$z_0.ah=ah$, for each $a\in A(S)$. Also each element of $A(S)^*$ is of the 
form $z_0T$, for some $T\in W^*(S)$, and one can then find a 
net $\{g_\alpha\}\subseteq C^*(S)$ such that $g_\alpha\overset{w^*}{\lo} T$ 
(since $W^*(S)$ is the envelping $W^*$-algebra of $C^*(S)$)
 to get $z_0g_\alpha\overset{w^*}{\lo}F$ and so 
\begin{align*}
<F,a\lambda> &=\lim_\alpha<z_0g_\alpha, a\lambda>=\lim_\alpha<g_\alpha, 
z_0.a\lambda>=\lim_\alpha<g_\alpha,a\lambda>\\
&=\lim_\alpha<g_\alpha, ah>=\lim_\alpha<g_\alpha, 
z_0.ah>=\lim_\alpha<z_0g_\alpha,ah>\\
&=<F,ah>.
\end{align*}
Therefore $\lambda=h$ as elements of $\tilde{M}(A(S))$. In particular $h\in 
B(S)$ and so $B(S)=\tilde{M}(A(S))=RM(A(S))=M(A(S))$. $\Box$

\end{document}